\documentclass{amsart}
     \usepackage{a4wide}
\usepackage{mathrsfs}
\usepackage{amssymb}                    

 \date{2000-07-30; 2001-01-12; 2001-09-24; 2002-08-08; 2004-07-28}

\newcommand{\C}{\mathscr{C}}
\newcommand{\D}{\mathscr{D}}
\newcommand{\F}{\mathscr{F}}
\newcommand{\N}{{\mathbb N}}

\newcommand{\Z}{{\mathbb Z}}
\newcommand{\on}{\restriction}

\newcommand{\Oclones}{{\mathscr O}}
\renewcommand{\P}{{\mathscr P}}

\def\xx#1 {\newtheorem{#1}[thm]{#1}}
\theoremstyle{remark}

\xx Acknowledgment
\xx Corollary
\xx Consequence
\xx Conclusion
\xx Construction
\xx Definition
\xx {Definition and Fact}
\xx Discussion
\xx Example
\xx {First Step}
\xx Idea
\xx Fact
\xx Lemma

\xx Notation
\xx {More Notation}
\xx Remark 
\xx Setup
\xx Suggestion
\xx Theorem

\xx {Main Lemma}

\xx{Open Questions}

\newcommand{\ee}{^{(1)}} 
\newcommand{\ff}{f\ee}

\newcommand{\nd}{\mathord{\raise1pt\hbox{$\nabla$}\!\!\Delta}}

\renewcommand{\P}{{\mathscr P}}

\makeatletter
\newcommand{\KNUTHcases}[1]
{\left \{\,\vcenter {\normalbaselines \m@th
\ialign {$##\hfil $&\quad ##\hfil \crcr #1\crcr }}\right .
}  
\makeatother

\def \itm#1 {\item[(#1)]}

\newcommand{\oo}[1]{\Oclones^{(#1)}}
\newcommand{\oc}[1]{\C^{(#1)}}
\newcommand{\od}[1]{\D^{(#1)}}

\newcommand{\E}{{\mathscr E}}
\renewcommand{\SS}{{\mathfrak S}}
\newcommand{\M}{{\mathscr M}}

\newcommand{\Pol}{{\rm Pol}}

\title{Large Intervals in the Clone Lattice}

\keywords{precomplete clones; local clones; dually atomic}

\author{Martin Goldstern}

\address{DMG/Algebra, TU Wien
\endgraf
Wiedner Hauptstra\3e 8-10/104.1
\endgraf
A-1040 Wien}
\email{Martin.Goldstern@tuwien.ac.at}
\urladdr{http://info.tuwien.ac.at/goldstern/}

\thanks{The first author is grateful to  the Hebrew University
of Jerusalem for the hospitality during his visit,
and  to the Austrian Science foundation for supporting the joint research
under FWF grant P13325-MAT}

\author{Saharon Shelah}
\thanks{The second author is supported by the
   German-Israeli Foundation for Scientific Research \& Development
   Grant No. G-294.081.06/93.   Publication number 747}

\address{Mathematics\endgraf Hebrew University of Jerusalem\endgraf 91904
Jerusalem,  Israel} 
\email{shelah@math.huji.ac.il}
\urladdr{http://math.rutgers.edu/\char`\~shelah}

\begin{document}

 \begin{abstract} 
 We give three examples of large intervals
 in the lattice of (local)
 clones on an infinite set $X$, by exhibiting clones $\C_1$, $\C_2$,
 $\C_3$ such that: 
 \begin{enumerate}
 \item  the interval $[\C_1,\Oclones]$ in the lattice of local clones is
 (as a lattice)
 isomorphic to $\{0,1,2,\ldots\}$ under the divisibility relation,
 \item the interval $[\C_2, \Oclones]$ in the lattice of local clones is 
 isomorphic to the congruence lattice of an arbitrary 
 semilattice,
\item the interval $[\C_3,\Oclones]$ in the lattice of all clones is 
 isomorphic to the lattice of all filters  on $X$. 
 \end{enumerate}
 These examples explain the difficulty of obtaining a satisfactory
 analysis of the clone lattice on infinite sets.  
 In particular, (1)  shows that the lattice of local clones is not
 dually atomic. 
 
\end{abstract}
\maketitle

\setcounter{section}{-1}

 \section{Introduction} \label{section.intro}

\begin{Definition}  
Let $X$ be a nonempty set.   The {\em full clone} on~$X$, called $\Oclones$
or $\Oclones(X)$ is the set of all finitary 
functions or (``operations'') on~$X$:  $\Oclones = \bigcup_{n=1}^ \infty \oo n$, 
where $\oo n$ is the set of
all functions from $X^n$ into~$X$. 
\\
A {\em clone} (on $X$)
 is a set $\C \subseteq \Oclones $
  which contains all
projections and is closed under composition. 


Alternatively, $\C$ is a clone if $\C$ is the set of term functions of
some universal  algebra over~$X$. 

For any set $\C \subseteq \Oclones$, we write $cl(\C)$ for the smallest
clone containing $\C$.  

\end{Definition}

The set of clones over $X$ forms a complete algebraic lattice with
largest element~$\Oclones$.  The coatoms of this lattice are called
``precomplete  clones'' or ``maximal clones''.  (See also \cite{Szendrei:1986},
\cite{PK:1979}).

\begin{Definition} A clone $\C$ is called a {\em local clone}, iff
each set $\C\cap \oo  k$ is closed in the product topology (Tychonoff
 topology) on
$X^{X^k}$, where $X$ is taken to be discrete. In other words, $\C$ is
local iff:
\begin{quote}
Whenever $f\in \oo k$, and for all finite sets $A \subseteq  X^k$ there is
$g\in \C$ with $g\on A= f\on A$, then $f\in \C$
\end{quote}
For any $\C \subseteq \Oclones$, we write $loc(\C)$ for the smallest local
clone containing $\C$.  (So $cl(\C) \subseteq loc(\C)$.) 
\end{Definition}

The set of local clones over $X$ forms again 
a complete lattice with
largest element~$\Oclones$ (a sublattice of the lattice of all clone).
  The coatoms of this lattice are called
``precomplete local clones'' or ``maximal local clones''.   Note that
not every precomplete local clone is also a precomplete clone.

If $X$ is finite (so  the notions of ``clone'' and ``local clone''
coincide), the precomplete clones are useful as a completeness
criterion. By a theorem of Rosenberg, 
\begin{itemize}
\item there are only finitely many (how many depends on the size of $X$)
precomplete clones, and in fact
there is an explicit list $\M_1^X, \ldots, \M_m^X$ of them, together with
 effective  procedures for checking $f\in \M_i$
\item every clone $\C \not=\Oclones$ is contained in a precomplete clone. 
\end{itemize}
Hence, there is an effective procedure 
to check if a given set of functions generates all of $\Oclones$,
i.e., $cl(\{f_1, \ldots, f_n\}) = \Oclones$:  Just check if  $\{f_1,\ldots ,
f_n\}$ is contained in one of the maximal clones.

For an infinite set $X$, it is known that there are infinitely many
precomplete clones (in fact: $2^{|X|}$ many precomplete local clones,
and $2^{2^{|X|}}$ many precomplete clones).   The question whether 
every clone  $\not=\Oclones$ on an infinite set $X$ is below a 
precomplete clone [``Is the lattice of clones dually atomic?'']
is  still not fully resolved (the answer is ``no'', if $X$ 
is countable and the continuum hypothesis holds, see \cite{GoSh:808}).

We will show here that the corresponding question for {\em local}
clones has a negative answer in ZFC. 

\cite{Rosenberg+Szabo:1984} showed that
there are  unbounded
chains in the lattice of local clones (i.e., chains whose only upper
bound is the trivial bound $\Oclones$). 

However, note that for a partial order $(P,{\le})$  the properties 
\begin{itemize}
\item [$(*)$] Every element of $P$ is below a maximal element
\item [$(**)$] Every chain of $P$ is bounded  
\end{itemize}
are in general not equivalent, not even if we restrict our attention
to those  partial orders which are of the form $L\setminus \{\max L\}
$, for a complete algebraic lattice $L$. 

The property $(**)$ always implies $(*)$ (by Zorn's lemma), and the
property $(*)$ trivially implies $(**)$ {\em if $P$ has only finitely
many maximal elements}.

\section{The lattice of local clones is not dually atomic}

On any infinite set $X$ we will define a local clone $\C_1$ such that
the interval $[\C_1,\Oclones]$ in the local clone lattice is isomorphic to
the natural numbers ordered by the divisibility relation.

%

%

\begin{Setup}
Fix an infinite set $X$, and let 
 $s:X \to X$ be  a 1-1 onto map without cycles.   In other words, $X =
 Y \times \Z$, and $s(y,n) = (y,n+1)$ for all $y\in Y$, $n\in \Z$. 
The orbits of $s$ (or: the sets $\{y\}\times \Z$ are called
 ``components''. 
\end{Setup}

\begin{Notation} We will write $\bar a$ or $\bar b$ to denote tuples 
$(a_1, \ldots, a_n)$ or $(b_1, \ldots b_k)$ (the values of $n$ or $k$
will be either irrelevant, or clear from the context).   

If $\bar a$, $\bar b$ are as above, then $(\bar a, \bar b)$ or 
$(\bar a, b_i: 1\le i \le k)$ denotes the $n+k$-tuple $(a_1, \ldots,
a_n, b_1, \ldots, b_k)$. 

For $n>0$, $s^n$ is the $n$-th iterate of $s$, $s^{-n}$ is the inverse
of $s^n$.  $s^0$ is the identity function. 
For $n\in Z$, $a\in X$ we may write $a+n$ instead of
$s^n(a)$.

More generally,  if $\bar a = (a_1,\ldots, a_k)$, $n\in \Z$, then we write
$\bar a + n$ for $  (s^n(a_1),\ldots, s^n(a_k))$.   Similarly, we
may write $ \bar a - n$ for $\bar a + (-n)$.

Clearly, $(\bar a + n_1) + n_2 = \bar a + (n_1 + n_2)$, and 
$(\bar a + n_1) - n_2 = \bar a + (n_1 - n_2)$,
 so we will
often omit parentheses. 
\end{Notation}

\begin{Definition and Fact} For every $n\in \N$, the set 
$$\Pol(s^n):= \bigcup_{k=1}^\infty \{f\in \oo k: \forall
\bar a [ f(\bar a + n ) = f(\bar a ) + n ] \} $$ is a local clone. 

For $n=0$, $\Pol(s^0)  = \Oclones$. 
\end{Definition and Fact}

These clones, and also the unbounded  chain $\Pol(s^{2^n})$ were already considered in \cite{Rosenberg+Szabo:1984}.

\begin{Theorem} \label{theorem.s} Let $s$ be as above. $\C_1:= \Pol(s)$. 
 Then the map $n\mapsto Pol(s^n)$
is a lattice isomorphism between the following two  lattices: 
\begin{itemize}
\item   $(\N, {|})$, the natural numbers with the divisibility
relation, where  $1$ is the smallest and $0$ the greatest element
\item $([\C_1, \Oclones],{\subseteq})$,  the set of all local clones extending
$\C_1$; this set is an
 interval in the lattice of local clones on $X$. 
\end{itemize}
In particular, there is no precomplete  local clone above $\C_1$.  Also, if $s$
has infinitely many components, then $(X, s) \simeq (X, s^n)$  for 
all $n\not= 0$, 
so all clones in $[\C_1, \Oclones)$ are
isomorphic (i.e., conjugate to each other via permutations of $X$). 
\end{Theorem}

\begin{Definition} 
We say that $\bar a$ and $\bar b$ are parallel (or $s$-parallel), 
  $$ \bar a \Vert \bar b \qquad  (\  \bar a \Vert_s \bar b \  ) 
$$
iff there is some $n\in \Z$, $\bar a + n = \bar b$.  

\end{Definition}

\begin{Fact}\label{fact.rep}
\begin{enumerate}
\item For each $k$, $\Vert$ is an equivalence relation on $X^k$. 
\item If $A \subseteq X^k$ meets each $\Vert$-equivalence
 class in at most one element [exactly one element], and $g:A \to X$, 
then there is a
 function [there is a unique function] $f\in \Pol(s)^{(k)}$
 with $f\on A = g$.
 \end{enumerate}
\end{Fact}

\begin{Definition}\label{def.a.n}
For any local clone  $\D \supseteq \Pol(s)$ define $G_\D$ and 
$n_\D$ as follows: 
\begin{align*}
 G_\D & := \{ n\in \Z:  \od 1 \subseteq \Pol(s^n) \}\\
 n_\D &:= \min \{ n\in G_\D, n> 0 \}, \qquad\qquad n_\D:=0 \mbox{  if
 $G_\D=\{0\}$} 
\end{align*}
\end{Definition}

\begin{Fact} $G_\D$ is a subgroup of $\Z$, hence $G_\D = \{ n_\D \cdot
 k: k\in \Z \}$
\end{Fact}
\begin{proof}
Let $n,m\in G_\D$, $f\in \od 1$.   We have to check 
$f(x+n-m) = f(x) + n-m$ for all $x\in X$:
$$f(x+(n-m)) = f((x + n)-m) = f(x+n)-m = f(x)+n-m$$
\end{proof}

\begin{Lemma}\label{lemma.1}
If $\D \supseteq \Pol(s)$ is a clone, and $\D\not= \Pol(s)$, then
already $\od 1 \not= \Pol(s)^{(1)}$. 
\end{Lemma}
\begin{proof}
Let $g\in \D\cap \oo k \setminus \Pol(s)$,
say $g(\bar a + 1 ) \not= g(\bar a) + 1$. Fix any $b\in X$.
   We can find unary functions
$f_1, \ldots, f_k\in \Pol(s)$ such that $f_i(b) = a_i$. 
\\Now consider the function $h:X\to X$, defined by 
$$ h(x) = g(f_1(x), \ldots, f_k(x)) $$
Clearly $h \in \od 1$, and $h(b+1) = g(\bar a + 1 ) \not= 
g(\bar a ) + 1 = h(b) + 1$, so $h\notin \Pol(s)$. 
\end{proof}

\begin{Lemma}
Let $\E \subseteq \oo 1 $,  and let $\D$ be the local clone generated
by $\Pol(s) \cup \E$, 
 $n^*:= n_\D$. 

If $\bar a = (a_1, \ldots, a_k)$ and $\bar b = (b_1, \ldots, b_k)$ are
not $s^{n^*}$-parallel, then there is a function $f\in \D$ such that 
$(\bar a, f(\bar a))$ and $(\bar b, f(\bar b))$ are not $s$-parallel. 
\end{Lemma}

\begin{proof}
First, note that $\Pol(s) \cup \od 1 \subseteq \Pol(s^{n^*})$,
so also $\D \subseteq \Pol(s^{n^*})$. 

If $\bar a$ and $\bar b$ are not $s$-parallel, then there is nothing
to do, so assume $\bar b = \bar a + \ell$. By our assumption, 
$\ell$ is not divisible by $n^*$, so $\ell\notin G_\D$. 

We can find a function $d\in \od 1$ and some $c\in X$ with
$d(c+\ell ) \not= d(c) + \ell$. \\
By fact \ref{fact.rep} there is a function $g\in \Pol(s)$ with
$g(a_1) = c$ (and $g(a_1 + \ell) = c+ \ell$). 

Now let $f:= d\circ g \circ \pi^n_1$, then $f(\bar a ) = d(c)$, 
$f(\bar a + \ell) = d(c+\ell) \not= d(c)+\ell = f(\bar a) + \ell$. 
\end{proof}

\begin{proof}[Proof of theorem \ref{theorem.s}]
It is clear that the map  $n \mapsto \Pol(s^n)$ maps natural numbers
to local clones above $\Pol(s)$, and that this map is 1-1.   It
remains to show that this map is onto. 

So let $\E \supseteq \Pol(s)$ be a local clone.  We will first consider
the clone $\D = loc(\E^{(1)} \cap \Pol(s))$ and prove that $\D=
\Pol(s^n)$ for some $n\in \N$.   If $n=0$, then $\D=\E=\Oclones$, and if
$n>0$ then we invoke lemma~\ref{lemma.1} to show $\E=\D$.

So we are now looking at a clone $\D\supsetneq \Pol(s)$, where $\D 
= loc(\od 1 \cup \Pol(s))$. 

Define $G_\D$ and $n^*:= n_\D$ as in definition \ref{def.a.n}.

We will prove $\Pol(s^{n^*} ) \subseteq \D$.
 Since $\D$ is a local clone, it is enough
to show that every function in $\Pol(s^{n^*})$ can be interpolated by
 a function in $\D$ on any finite set. 

 So let $g\in
\Pol(s^{n^*})$ be $k$-ary. 

 Let $\bar a_1$, \dots, $\bar a_n$ be in
$X^k$, $b_\ell:= g(\bar a_\ell)$.   We claim that there is a function
$f\in \D$ satisfying also 
 $b_\ell:= f(\bar a_\ell)$ for $\ell= 1,\ldots,n$. 

Wlog we may  assume that no two of the $k$-tuples $\bar a_\ell$ are
$s^{n^*}$-parallel.  [This assumption is allowed, since $g \in
\Pol(s^{n^*})$.]

Let $I$ be the set of all pairs $i = (\ell_1, \ell_2) $ of distinct
numbers in $\{1, \ldots, n\}$.   For each $i= (\ell_1, \ell_2)\in I$ 
we can find a function $f_i\in \D$ such that 
$(\bar a_{\ell_1}, f_i(\bar a_{\ell_1}))$ and 
$(\bar a_{\ell_2}, f_i(\bar a_{\ell_2}))$ are not $s$-parallel.

Let $\bar c_\ell:= (\bar a_\ell, f_i(\bar a_\ell): i \in I)$. Clearly, 
for all $\ell_1 \not= \ell_2$ we have: \\
$\bar c_{\ell_1}$ and  $\bar c_{\ell_2}$ are not $s$-parallel. 

So by fact \ref{fact.rep} there is a function $h\in \Pol(s)$ such that 
$$ \forall \ell: \ \ h(\bar c_\ell) = b_\ell $$
So the function $f$ defined by  
$$f(\bar x) = h(\bar x, f_i(\bar x): i \in I)$$
is in $\D$, and it satisfies $f(\bar a_\ell)  = b_\ell$ for all
$\ell$. 

\end{proof}

\begin{Remark}   Let $H$ be any constant function.  
Then  $cl(\Pol(s) \cup \{H\}) = \Oclones$.
    It is well known (and easy to show,
using Zorn's lemma) that this implies that every clone above
$\Pol(s)$ is below some precomplete clone. 

(In particular, also the (nonlocal) clone 
 $\C_\infty:= \bigcup_{n > 0 } \Pol(s^n)$ is below a precomplete
 clone.) 
\end{Remark}

\section{A large interval of local  clones}

\begin{Theorem}
Let $\SS = (S,\vee)$ be a downward directed semilattice,
 and let $Con(\SS)$ be
the lattice of congruences on $\SS$.  Then there is a local clone
$\C_2$ (on the set~$S$) such that 
$$ [\C_2, \Oclones_S ] \simeq Con(\SS)$$
That is, there is a lattice isomorphism between the set of local
clones above $\C$ and the set of congruences of $\SS$. 
\end{Theorem}

Remark: If $\emptyset \subsetneq I \subsetneq S$ is an ideal, then
the partition $\{ I, S \setminus I \}$ corresponds to a congruence
relation which is a coatom in $Con(\SS)$.  In fact, all coatoms are
obtained in this form.   It is clear that $Con(\SS)$ is dually
atomic.

\newcommand{\pres}{\Pol}
\begin{Definition} Let $S$ be a set $A \subseteq S$.  We let $\pres (A)$ be
the set of all functions on $S$  which ``preserve''~$A$: 
$$ \pres (A) := \bigcup_{n=1}^\infty \{f\in \oo n : f[A^n] \subseteq A\}$$
\end{Definition}

\begin{Fact}
$\pres (A)$ is always a local clone. 

If $\emptyset \not= A \not= S$, then $\pres (A)$ is a maximal clone. 
\end{Fact}

\begin{Definition and Fact}\label{def.corder}
Let $\SS=(S,\vee)$ be a semilattice.  We call $R \subseteq S\times S$ 
a {\em congruence order} on $\SS$ 
iff one (or both) of the following two equivalent 
conditions are satisfied: 
\begin{enumerate}
\item ${\theta} _R:= \{(x,y): x R y \mbox { and } y R x \}$ is a
(semilattice) congruence relation, and: $x R y  $ iff $x/{\theta}_R
\le y / {\theta}_R$. 
\item $R$ is reflexive and transitive, $x\le y \Rightarrow x R y$, and
$$(*) \qquad \qquad \forall x,y,z: \ \ 
 x R z \,\, \& \,\, y R z \ \Rightarrow \ (x \vee y) R z $$
\end{enumerate}
\end{Definition and Fact}
The following fact is trivial: 
\begin{Fact} The maps $ R \mapsto {\theta}_R$ and ${\theta}  \mapsto 
\{(x,y): x/{\theta}  \le y/{\theta} \}$ are monotone 
bijections between congruence
relations and congruence orders, and they are inverses of each other. 
\end{Fact}

\begin{Notation}
For $a,b\in S$ let $\chi_{a,b}$ be the function satisfying 
$\chi_{a,b}(b)=a$, $ \chi_{a,b}(x)=x$ for $x\not= b$. 
\end{Notation}

\begin{Notation} For $\bar a = (a_1,\ldots, a_k)$, write $\bigvee 
\bar a $ for $a_1 \vee \cdots \vee a_k$. 
\end{Notation}

\begin{Definition and Fact} \label{def.sq} 
Let $\SS = (S,{\vee})$ be a semilattice, 
and let  $\sqsubseteq $ be a congruence order
on $\SS$. Then
$$ \E({\sqsubseteq}):=  
\bigcup_{k=1}^\infty
\{ f\in \oo k: \forall \bar x \in S^k\ \bigl[f(\bar x)
\sqsubseteq \bigvee \bar x \bigr] \}$$ is a local clone. 
Furthermore,  
\begin{enumerate}
\item
$ \E({\sqsubseteq}) = \bigcap_{a\in S}  \pres{\{ x\in S: x\sqsubseteq
a \}} $
\item ${\sqsubseteq_1} \subseteq { \sqsubseteq_2}$ implies 
        $\E({\sqsubseteq_1}) \subseteq \E({\sqsubseteq_2})$. 
\item   $ f\in \E(\sqsubseteq )\cap \oo 1 $ implies $\forall x:
f(x)\sqsubseteq x$. 
\item $\chi_{a,b}\in \E(\sqsubseteq)$ iff $a \sqsubseteq b$. 
\end{enumerate}
\end{Definition and Fact}

We will now consider local clones above the clone $ \C_2 := \E(\le)$ and
we will show that they all are induced by congruence orders/congruence
relations, and that also conversely every congruence relation is
induced by a clone. 

\begin{Definition and Fact}
Let $\C \supseteq \C_2$ be a local clone. Then 
$$ {\sqsubseteq_\C } := \{(x,y):  \exists f\in \oc 1\, f(y)=x\}$$
is a congruence order. 
\\
Also, $\C_2 \subseteq \C_1 \subseteq \C_2$ implies 
${ \sqsubseteq_{\C_1}} \subseteq  {\sqsubseteq_{\C_2}} $. 
\end{Definition and Fact}
\begin{proof}  We will check condition (2) from definition
\ref{def.corder}.  Clearly $\sqsubseteq$ is reflexive and transitive. 
If $a\le b$, then the function $\chi_{a,b}\in \C_2 \subseteq \C$
will witness  that $a \sqsubseteq_\C b$. 

It remains to check $(*)$.   So let $a\sqsubseteq c$, $b\sqsubseteq
c$. There are functions $f,g\in \C$ with $f(c)=a$, $g(c)=b$. 
Since the function $\vee: (x,y) \mapsto x\vee y$ is in $\C_2 \subseteq
\C$, we also have $(f\vee g)\in \C$, and $f\vee g$ witnesses $a \vee b
\sqsubseteq_\C  c$. 
\end{proof}

\begin{Lemma} Let $R$ be a congruence order.  Then $R =
{\sqsubseteq_{\E(R)}}$. 
\end{Lemma}
\begin{proof} 
The inclusion ${\sqsubseteq_{\E(R)}} \subseteq R$ is
trivial:  Let $a \sqsubseteq_{\E(R)} b$.  So there is $f\in \E(R)$,
$f(b) = a$. Now  by fact \ref{def.sq}(3), $b R a$.

For the proof of the reverse inclusion, 
 $R \subseteq {\sqsubseteq_{\E(R)}}$, consider any 
 $ a R b$. The function $\chi_{a,b}\in \E(R)$
 witnesses 
 $a \sqsubseteq_{\E(R)} b$. 
\end{proof}

\begin{Lemma} Let $\C \supseteq\C_2$ be a local clone.  Then $\C =
\E(\sqsubseteq_\C)$. 
\end{Lemma}

\begin{proof}
The inclusion $\C \subseteq 
\E(\sqsubseteq_\C)$ is trivial:   For  $f\in \C$ we need to show that
for all  $a\in S$, $f$ preserves the set $\{x: x \sqsubseteq_\C a\}$. 
Let $\bar x = (x_1,\ldots, x_k)$. 
If $x_1, \ldots, x_k \sqsubseteq_\C a$
 then for each $i$ there  is some $g_i\in\C$ with $g_i(a)=x_i$.  Now 
$f(g_1,\ldots, g_k)$ witnesses that also $f(\bar x)\sqsubseteq_\C a$. 

Now we will show 
$\E(\sqsubseteq_\C) \subseteq \C$: \\
Let $f\in \E(\sqsubseteq_\C)$ be $k$-ary,
 where $\C \supseteq \C_2$. To show that $f\in
\C$ it is enough [since $\C$ is local] 
to show that $f$ can be interpolated by an element of
$\C$ on any finite number of places . 

So let $\bar a_1,\ldots, \bar a_n\in S^k$, and let $b_i:= f(\bar
a_i)$.    Let $d_i:= \bigvee \bar a_i$.     

Since $f\in \E(\sqsubseteq_\C)$, we have $b_i \sqsubseteq_\C d_i$, so there is (for
every $i$) a
unary function $g_i\in \C$ with $g_i(d_i)= b_i$. 

Define a $k+1$-ary function $h_i$ by letting $h_i(y, \bar x) = y$ if
$\bar x = \bar a_i$, and 
\begin{quote}
if $\bar x \not= \bar a_i$, then: \\ 
$h(y, \bar x):= $ some value which is $\le x_j$ and $\le b_j$
for all $j$, and also $\le y$.  
\end{quote}
(It is  possible  to find
such a value, since  $\SS$
is downward directed.) \\
 Clearly $h_i\in \C \subseteq
\C$. So the function 
$$ f_i: \bar x \mapsto h_i(g_i(\bigvee \bar x), \bar x)$$
is in $\C$.   Now check that the function $f':= \bigvee\limits_{i=1}^n
f_i$ maps $\bar a_i$ to $b_i$. Clearly  $f'\in \C$. 
\end{proof}

\begin{Example}
Let $(S,{<})$ be a linearly ordered set.  Then the congruence relations on
$(S,\max)$ are exactly the equivalnce relations with convex classes. 
\end{Example}

\begin{Example}
As a special case, consider the semilattice $(\N,\max)$.  A congruence relation is just a partition of $\N$ into disjoint intervals. 

The map $$\theta \mapsto A_\theta:= \{
\max E:  \hbox{$E$ is a finite congruence class}\}$$
is an antitone 1-1 map from the congruence relations into $\P(\N)$, the 
power set of $\N$. 

   It is also easy to see that this map is onto: Each $A \subseteq \N$ is 
equal to $A_{\theta(A)}$, where for $k<n$ we have: 
\begin{quote}
$ (k,n)\in \theta(A) $ iff there is no $a\in A$, $k\le a< n$
\end{quote}

The map $A\mapsto \bigcap_{a\in A} \Pol\{0,\ldots, a\}$ is an
isomorphism between $(\P(\N), {\supseteq})$ and $[\C_2, \Oclones]$.  
The empty set corresponds to $\Oclones$, or to the equivalence relation 
with a single class;  the set $\N$ itself corresponds to $\C_2$, or to the 
 equivalence with singleton classes. 

\end{Example}

%
%
%
%
%

\section{A large interval of clones}

\def\fix{\operatorname{fix}}
\def\nix{\operatorname{nix}}

On any infinite set $X$ we will define a clone $\C_3$ such that the
 interval $[\C_3,\Oclones]$ in the full clone lattice is very large (with
 $2^{2^{|X|}}$ precomplete elements), but still reasonably well understood.


\begin{Definition}  For any function $f\in \oo n$, let $\ff\in \oo 1 $
be defined by $\ff(x) = f(x,\ldots, x)$. 
\end{Definition}
\begin{Definition}  For any function $f\in \oo n$, we let 
$$ \fix(f) = \{x: \ff(x) = x\} \qquad 
 \qquad 
\nix(f) =  \{x: \ff(x) \not= x\}$$
\end{Definition}

\begin{Definition}  Let $\F \subseteq \P(X)$ be a family of sets. 
We define 
$$\C_\F:= \{ f\in \Oclones:  \fix(f) \in \F \}$$
\end{Definition}

\begin{Fact}
\begin{enumerate}
\item   If $\F$ is a filter, then $\C_\F$ is a clone. 
\item  If $\F \subseteq \F'$, then $\C_\F \subseteq \C_{\F' }$. 
\end{enumerate}
\end{Fact}

\begin{Definition}  Let $\C_3 := \C_{\{X\}}$ be the clone of 
``idempotent'' functions, 
i.e., of all functions satisfying  $f(x,\ldots, x) = x$ for all $x$. 
\end{Definition}

\begin{Theorem}  
The map $\F \to \C_\F$ is an order isomorphism between the set of all
filters  (including the improper filter $\P(X)$)
and the set of all clones above $\C_3$. 

In particular, the precomplete clones above $\C_3$ are exactly the
clones of the form $\C_U$, where $U$ is an ultrafilter on $X$. 
\end{Theorem}

We will prove this theorem in several steps, concluding with 
 lemma \ref{lemma.final} below.

\begin{Lemma}  Let $\F$ be a filter, $\D \supsetneq \C_\F$.   Then 
$\D^ {(1)}\supsetneq  \C_\F^{(1)}$. 
\end{Lemma}
\begin{proof}
Let $f\in \D \setminus \C_\F$.  Then $\fix(f)= \fix(\ff) \notin \F$,
so $\ff\in \D^{(1)} \setminus \C_\F$. 
\end{proof}


\begin{Lemma}\label{main} Assume $f\in \D\supseteq \C_3$, and $\fix(f) \subseteq
\fix(g)$,  $f,g\in \oo
1$.   Then $g\in \D$. 
\end{Lemma}

\begin{proof}

Let $$H(x,y) = \KNUTHcases{ g(x) & if $x\not=y$\cr
                            x    & if $x=y$\cr}$$
Clearly $H\in \C_3$.  For $x\in \fix(f) \subseteq \fix(g)$ we have
           $H(x,f(x))=H(x,x)=x=g(x)$, and 
for $x\in \nix(f)$ we have $f(x)\not=x$, so $H(x,f(x))=g(x)$. 

So in either case, $H(x,f(x))=g(x)$. 
\end{proof}

\begin{Lemma} \label{lemma.final}
Let $\D$ be a clone with $\C_3\subseteq \D$, $\D\not=
\Oclones$. 
  Then there is a (proper)
filter $\F$ such
that $\D= \C_\F$. 
\end{Lemma}
\begin{proof}
Let $I:= \{\nix(f): f\in \D\} = \{\nix(f): f\in \D^{(1)}\}$.   We
first check that $I$ is an ideal. 

If $A= \nix (f)$, $f\in \D$, and $B \subseteq A$, then by lemma
\ref{main} there is a
function $g\in \D$ 
with $B= \nix(g)$.  So $I$ is downward closed. 

Now let $A_\ell = \nix(f_\ell) $, $f_\ell\in \D$
 for $\ell=1,2$, and assume that
$A_1\cap A_2= \emptyset$.   

Let $B = X \setminus (A_1\cup A_2)$.   We may assume that either
$|A_1|\ge 2$, or $B\not=\emptyset$ (or both). 

In either case there is a unary function $f_1' $ with 
 $\nix (f_1') = A_1$, and $f_1'$ maps
$A_1$ into $A_1\cup B$.  
By lemma \ref{main}, $f_1'\in \D$. 
 So 
$f_2\circ f_1'\in \D$.  Also, $\nix(f_2\circ f_1') = A_1\cup A_2$;
this shows that $I$ is an ideal.

Let $\F$ be the filter dual to $I$.  Clearly, 
$$ f\in \D \Rightarrow \fix f \in \F \Rightarrow f\in \C_\F. $$ 

For the converse, we first check $\C_\F\cap \oo 1  \subseteq \D$: 
\\
Let $f\in \C_\F$ be unary.  So $\fix( f)\in \F$, i.e., there is a
function $g\in \D$ with $\fix f = \fix g$.  By lemma \ref{main}, $f\in
\D$. 

Now take an arbitrary $n$-ary function $f\in \C_\F$. We need to show
that $f\in \D$.     Let $A= \nix f\in I$. 
Define an $n+1$-ary function $H$ as follows: 
$$H(x_1,\ldots, x_n, y) = \KNUTHcases{ x_1 & if $x_1=\cdots = x_n=y$\cr
f(x_1, \ldots x_n) & otherwise}$$
Clearly $H\in \C_3$.  Note that $f\ee\in \D$, so also the function 
 $\bar x \mapsto H( \bar x , f\ee(x_1)) $ is in $\D$.

We now check that $H(\bar x, f\ee(x_1)) = f(\bar x)$ for all $\bar x$. 
We distinguish three cases: 

\begin{itemize}
\item[Case 1:] $x_1=\cdots = x_n \in A=\nix(f)$.  So $f\ee(x_1) \not=
x_1$, hence (by definition of $H$) we have
$ H(\bar x, f\ee(x_1)) = f(\bar x)$. 
\item[Case 2:]  $x_1=\cdots = x_n \in \fix(f)$.  So $f(\bar x) =
f\ee(x_1) =x_1$, and also $H(\bar x, f\ee(x_1)) = x_1$. 

\item[Case 3:] Not all $x_i$ are equal. Again, by definition of $H$, 
we have  $H(x_1, \ldots, x_n, f\ee(x_1)) = f(x_1,\ldots, x_n)$.  
 \end{itemize}

This shows that $f\in \D$. 

\begin{Remark}
If we regard the set $X$ as a discrete topological space, then the
Stone-Cech compactification of $X$ is 
$$ \beta X = \{ U:  \mbox{$U$ is an ultrafilter on $X$}\}$$
There is a canonical 1-1 order-preserving correspondence between the
filters on $X$ (ordered by $\subseteq$)  and 
the closed subsets of $\beta X$ (ordered by $\supseteq$).

So the interval $[\C_3,\Oclones]$ in the full clone lattice  is isomorphic (as a
complete lattice)  to the
family of closed subsets of $\beta X$, ordered by reverse inclusion:
$\Oclones$ corresponds to the empty set, each precomplete clone in $[\C_3,\Oclones]
$ corresponds to a singleton set.

Note that for any closed subset $F \subseteq \beta X$ and any
$p\in \beta X \setminus F$, also $F \cup \{p\}$ is closed, and 
moreover: 
\begin{quote}
 $F $ covers $G$ (i.e., $F \supset G$, and the interval $(G,F)$ is 
empty) iff
$G = F\cup \{p\}$ for some $p\in \beta X \setminus F$
\end{quote}

In particular, let $\C_{\rm bd}\supseteq \C_3 $ be the clone corresponding to 
the ideal of small sets, i.e., 
$$\C_{\rm bd} := \{ f\in \oo : \exists B \subseteq X, |B| < |X|, \forall x\in X\setminus B: f(x,\ldots, x) = x\}$$

Then every clone $\C \supsetneq \C_{\rm bd}$ has exactly $2 ^{2^{|X|}}$ 
lower neighbors in the clone lattice; the clone corresponds to a closed
set $F$, and the lower neighbors correspond to closed sets $F\cup \{p\}$. 

This is a special case of a theorem of  \cite{Marchenkov:1981}. 
\end{Remark}

\end{proof}

\nocite{Freese+Nation:1973}
\nocite{Zitomirski:1971}

\bibliographystyle{lit-unsrt}



\end{document}